\documentclass[12pt]{amsart}
\usepackage{amsmath,amssymb,amsthm,amscd}
\newtheorem{thm}{Theorem}
\newtheorem*{thmnonum}{Theorem}

\newtheorem{cor}[thm]{Corollary}
\newtheorem*{rem}{Remark}

\newcommand{\SL}{{\rm SL}}
\newcommand{\C}{\mathbb{C}}
\newcommand{\Z}{\mathbb{Z}}
\newcommand{\R}{\mathbb{R}}
\newcommand{\Aut}{{\rm Aut}}
\newcommand{\IM}[1]{{\rm Im}~#1}
\newcommand{\abs}[1]{\left\vert#1\right\vert}

\begin{document}

\title{Bounds for coefficients of cusp forms and extremal lattices}
\author{Paul Jenkins}
\address{Department of Mathematics, Brigham Young University, Provo, UT 84602}
\email{jenkins@math.byu.edu}

\author{Jeremy Rouse}
\address{Department of Mathematics, Wake Forest University,
  Winston-Salem, NC 27109}
\email{rouseja@wfu.edu}
\thanks{The second author was supported by NSF grant DMS-0901090}
\begin{abstract}
A cusp form $f(z)$ of weight $k$ for $\SL_{2}(\Z)$
is determined uniquely by its first $\ell := \dim S_{k}$ Fourier
coefficients. We derive an explicit bound on the $n$th coefficient
of $f$ in terms of its first $\ell$ coefficients. We use this result to
study the non-negativity of the coefficients of the unique modular form of weight $k$ with
Fourier expansion
\[
  F_{k,0}(z) = 1 + O(q^{\ell + 1}).
\]
In particular, we show that $k = 81632$ is the largest weight for
which all the coefficients of $F_{k,0}(z)$ are non-negative. This
result has applications to the theory of extremal lattices.
\end{abstract}

\maketitle

\section{Introduction and Statement of Results}

An incredible number of interesting sequences appear as Fourier
coefficients of modular forms.  The analytic properties of these modular forms
dictate the asymptotic behavior of the corresponding sequences.

The most famous example of such a sequence is the partition function $p(n)$,
which counts the number
of ways of representing an integer $n$ as a
sum of a non-increasing sequence of positive integers. Hardy and
Ramanujan pioneered the use of the circle method to study the
asymptotics for $p(n)$ and proved that
\[
  p(n) \sim \frac{1}{4n \sqrt{3}} e^{\pi \sqrt{\frac{2n}{3}}}
\]
by using the analytic properties of the generating function
\[
  f(z) = \sum_{n=0}^{\infty} p(n) q^{n} =
  \prod_{n=1}^{\infty} \frac{1}{1 - q^{n}}, \] where $q = e^{2 \pi i z}$.
(See Chapter 5 of \cite{Andrews} for a proof as well as for an exact
formula for $p(n)$).

Another important example is given by the arithmetic of quadratic forms.
Let $Q$ be a positive-definite, integral, quadratic
form in $r$ variables, where $r$ is even, and let $r_Q(n)$ denote the
number of representations of the integer $n$ by $Q$. It is well-known that the
generating function
\[
  \theta_{Q}(z) = \sum_{n=0}^{\infty} r_{Q}(n) q^{n}
\]
is a holomorphic modular form of weight $\frac{r}{2}$ for some
congruence subgroup of $\SL_{2}(\Z)$ (see Chapter 10 of \cite{Iwa}
for details).

To determine which integers are represented by $Q$, it is necessary
to study the decomposition
\[
  \theta_{Q}(z) = E(z) + G(z)
\]
where $E(z)$ is an Eisenstein series and $G(z)$ is a cusp form, and
to determine explicit bounds on the coefficients of $E(z)$ and
$G(z)$.  If $r \geq 6$, formulas for the coefficients of Eisenstein series
show that the coefficients of $E(z)$ are of size $n^{\frac{r}{2} - 1}$,
and if we write
\[
  G(z) = \sum_{i=1}^{\ell} c_{i} g_{i}(d_{i} z)
\]
where the $g_{i}(z)$ are newforms, then Deligne's proof of the Weil
conjectures implies that the $n$th coefficient of $G(z)$ is bounded
by
\[
  \left(\sum_{i=1}^{\ell} |c_{i}|\right) d(n) n^{\frac{r - 2}{4}}.
\]

In \cite{BH}, Bhargava and Hanke prove that a positive-definite
quadratic form with integer coefficients represents every positive
integer if and only if it represents the integers from 1 up to 290;
in fact, it is only necessary for the form to represent 29 of these
numbers. To prove this, they study about $6000$ quadratic forms in
four variables, and the most time-consuming part of their
calculation comes from computing the constant
\[
  C(G) = \sum_{i=1}^{\ell} |c_{i}|.
\]
In this paper, we find bounds for this constant $C(G)$ for general cusp forms
$G$ of weight $k$ and full level.

If
\[
  \ell := \dim S_{k} = \begin{cases}
  \lfloor \frac{k}{12} \rfloor & \text{ if } k \not\equiv 2 \pmod{12}\\
  \lfloor \frac{k}{12} \rfloor - 1 & \text{ if } k \equiv 2 \pmod{12},
  \end{cases}
\]
then any cusp form $G(z) = \sum_{n=1}^{\infty} a(n) q^{n}$ is
determined uniquely by the coefficients $a(1)$, $a(2)$, $\ldots$,
$a(\ell)$. In fact, in \cite[Theorem 3]{BKO}, Bruinier, Kohnen and Ono
showed that the coefficients $a(n)$ of $G(z)$ may be explicitly
computed recursively from the first $\ell$ coefficients of $G$.
Specifically, $a(n)$ may be written as a polynomial with rational
coefficients in the coefficients $a(n-i)$, the weight $k$, and the
values of the $j$-function at points in the divisor of $G$.

Our first result is a bound on $\sum_{i=1}^{\ell} |c_{i}|$ (giving a bound on
$\abs{a(n)}$) in terms of the coefficients $a(1), a(2), \ldots, a(\ell)$.

\begin{thm}
\label{bound} Assume the notation above. Then
\[
  |a(n)| \leq \sqrt{\log(k)} \left(11 \cdot \sqrt{\sum_{m=1}^{\ell}
  \frac{|a(m)|^{2}}{m^{k-1}}}
  + \frac{e^{18.72} (41.41)^{k/2}}{k^{(k-1)/2}}
  \cdot \left|\sum_{m=1}^{\ell} a(m) e^{-7.288m} \right| \right) \cdot d(n) n^{\frac{k-1}{2}}.
\]
\end{thm}

We apply this result to the study of extremal lattices. An even,
unimodular lattice is a free $\Z$-module $\Lambda$ of rank $r$,
together with a quadratic form $Q : \Lambda \to \Z$ with the property that the
inner product
\[
  \langle \vec{x}, \vec{y} \rangle = Q(\vec{x} + \vec{y}) - Q(\vec{x}) - Q(\vec{y})
\]
is positive definite on $\R \otimes \Lambda$ and is an integer for all pairs $\vec{x}, \vec{y} \in \Lambda$; additionally, we require that $\langle \vec{x},
\vec{x} \rangle$ is even for all $\vec{x} \in \Lambda$, and that
the dual lattice
\[
  \Lambda^{\#} := \{ \vec{x} \in \R \otimes \Lambda : \langle \vec{x}, \vec{y} \rangle \in \Z \text{ for all } \vec{y} \in \Lambda \}
\]
is equal to $\Lambda$. For such a lattice, we must have
$r \equiv 0 \pmod{8}$, so the theta function
$\theta_{Q}$ is a modular form for $\SL_{2}(\Z)$ of weight $k \equiv 0 \pmod{4}$.

For example, if
\[
  Q = x_{1}^{2} + x_{2}^{2} + x_{3}^{2} + x_{4}^{2} + x_{5}^{2} + x_{6}^{2}
  + x_{7}^{2} + x_{8}^{2} - x_{1} x_{3} - x_{2} x_{4} - x_{3} x_{4} - x_{4} x_{5} - x_{5} x_{6} - x_{6} x_{7} - x_{7} x_{8},
\]
then $\Lambda$ is the $E_{8}$ lattice and
\[
  \theta_{Q}(z) = E_{4}(z) = 1 + 240 \sum_{n=1}^{\infty} \sigma_{3}(n) q^{n}.
\]

An even, self-dual lattice $\Lambda$ is called extremal if $r_{Q}(n)
= 0$ for $1 \leq n \leq \lfloor \frac{r}{24} \rfloor$. This means
that if $Q$ is the quadratic form corresponding to $\Lambda$, then
\[
  \theta_{Q}(z) = 1 + O(q^{\ell + 1}) \in M_{\frac{r}{2}}.
\]
An example is given by the famous Leech lattice $\Lambda_{24}$. It
is the unique extremal lattice of dimension 24, and
$\Aut(\Lambda_{24})$ is a perfect group whose quotient by $-1$ is
$Co_{1}$, the first sporadic finite simple group discovered by John
H. Conway.

Little is known about the set of dimensions in which extremal lattices exist,
and examples are known only in dimensions $\leq 88$. Cases where the
rank is a multiple of $24$ are particularly challenging, and Nebe~\cite{Nebe}
recently succeeded in constructing a $72$-dimensional extremal
lattice.

If $\Lambda$ is an extremal lattice of dimension $r$, then the
definition of $r_Q(n)$ implies that all the Fourier coefficients
of the modular form
\[
  \theta_{Q}(z) = \sum_{n=0}^{\infty} r_{Q}(n) q^{n} = 1 + O(q^{\ell + 1})
  \in M_{\frac{r}{2}}
\]
are non-negative. In \cite{MOS}, Mallows, Odlyzko, and Sloane use this
to show that extremal lattices fail to exist in large dimensions
(larger than about 164,000) by showing that the unique modular form of weight $k$
with Fourier expansion
\[
  F_{k,0}(z) = \sum_{n=0}^{\infty} a(n) q^{n} = 1 + O(q^{\ell + 1}),
\]
has $a(\ell + 2) < 0$ if
$k$ is large enough. (In \cite{S}, Siegel proved that $a(\ell + 1) >
0$ for all $k \equiv 0 \pmod{4}$).

As an application of Theorem~\ref{bound}, we give an explicit
estimate on the largest index negative coefficient of $F_{k,0}(z)$.

\begin{thm}
\label{extreme} Suppose that $k \equiv 0 \pmod{4}$, and $F_{k,0}(z)
\in M_{k}$ is the unique modular form of weight $k$ with
\[
  F_{k,0}(z) = 1 + O(q^{\ell + 1}) = \sum_{n=0}^{\infty} a(n) q^{n}.
\]
We have $a(n) > 0$ if
\[
  n \geq e^{58.366/(k-2)} (\ell^{3} \log(k))^{\frac{1}{k-2}} 1.0242382 \ell.
\]
\end{thm}

\begin{rem}
The result above is surprisingly strong. The factor preceding
$1.0242382 \ell$ tends to $1$ as $k \to \infty$, and since $a(n) =
0$ for $n \leq \ell$, the only region in which negative coefficients
could occur is (asymptotically)
\[
  \ell < n < 1.0242382 \ell.
\]
\end{rem}

We now use this bound to determine the largest weights $k$ in which
all the coefficients of $F_{k,0}(z)$ are non-negative. This depends
on $k \pmod{12}$, and so we have three cases.

\begin{cor}
\label{cor0} The largest weight $k$ for which all coefficients of
$F_{k,0}(z)$ are non-negative is \begin{align*} k = 81288 \, \,
&\text{if} \, \, k \equiv 0 \pmod{12}, \\ k = 81460 \, \, &\text{if}
\, \, k \equiv 4 \pmod{12}, \text{and} \\ k = 81632 \, \, &\text{if} \, \, k
\equiv 8 \pmod{12}.
\end{align*}
\end{cor}

%

\begin{rem}
As a consequence, the largest possible dimension of an extremal
lattice is $163264$.
\end{rem}

Our approach to proving our results is to study the basis of cusp
forms
\[
  F_{k,m}(z) = q^m + \sum_{n=\ell + 1}^{\infty} A_k(m,n) q^{n}
  \in S_{k}.
\]
Theorem 2 of \cite{DJ} gives a generating function for the forms
$F_{k,m}(z)$, and by integrating this generating function
we are able to isolate individual coefficients of these forms. Using
this method leads to a bound of the form
\[
  |A_k(m,n)| \leq c_{1} \cdot c_{2}^{\ell} e^{c_{3} m + c_{4} n}
\]
where $c_{1}, c_{2} > 0$, $c_{3} < 0$ and $0 < c_{4} < \sqrt{3}/2$.
Given that the coefficients of a cusp form of weight $k$ are bounded
by $O(d(n) n^{\frac{k-1}{2}})$, this bound is not useful by itself.
Next, we estimate the Petersson norm $\langle F_{k,m}, F_{k,m}
\rangle$ which is (essentially) the infinite sum
\[
  \sum_{n=1}^{\infty} \frac{|A_k(m,n)|^{2}}{n^{k-1}}
  \int_{2 \pi \sqrt{3} n}^{\infty} y^{k-2} e^{-y} \, dy.
\]
The exponential decay in the integral now cancels the exponential
growth from the bound on $|A_k(m,n)|$. Finally, we translate the
bound on $\langle F_{k,m}, F_{k,m} \rangle$ to a bound on the
constant $\sum_{i=1}^{\ell} |c_{i}|$ using methods similar to those
in \cite{Rou}.

An outline of the paper is as follows. In Section~\ref{prelim} we
review necessary background material about modular forms. In
Sections~\ref{proofofbound} and \ref{proofofextreme} we prove
Theorems~\ref{bound} and \ref{extreme}, respectively.  In
Section~\ref{computations}, we prove Corollary~\ref{cor0}.

\section{Preliminaries}
\label{prelim}

Let $M_{k}$ denote the $\C$-vector space of all holomorphic
modular forms of weight $k$ for $\SL_{2}(\Z)$, and let $S_{k}$ denote
the subspace of cusp forms.  For even $k \geq 4$, we have the classical Eisenstein series
\[
  E_{k}(z) = 1 - \frac{2k}{B_{k}} \sum_{n=1}^{\infty} \sigma_{k-1}(n) q^{n} \in M_{k},
\]
where $B_k$ is the $k$th Bernoulli number and $\sigma_{k-1}(n)$ is the
sum of the $k-1$st powers of the divisors of $n$.  We will also use the standard $\Delta$-function
\[\Delta(z) = \frac{E_4^3-E_6^2}{1728} = q \prod_{n=1}^\infty (1-q^n)^{24} = \sum_{n=1}^\infty \tau(n) q^n \in S_{12}\]
and the classical modular $j$-function
\[j(z) = \frac{E_4(z)^3}{\Delta(z)} = q^{-1} + 744 + 196884 q + \ldots,\]
a weakly holomorphic modular form of weight 0.  (Weakly holomorphic modular forms are holomorphic
on the upper half plane and satisfy the modular equation, but may have poles at the cusps.)

For each prime $p$, there is a Hecke
operator $T_{p} : M_{k} \to M_{k}$ given by
\[
  \sum_{n=1}^{\infty} a(n) q^{n} | T_{p}
  := \sum_{n=1}^{\infty} \left(a(pn) + p^{k-1} a\left(\frac{n}{p}\right)\right)
  q^{n}.
\]
The subspace $S_{k}$ is stable under the action of the Hecke operators.

If $f, g \in S_{k}$, we define the Petersson inner product of
$f$ and $g$ by
\[
  \langle f, g \rangle
  = \frac{3}{\pi} \int_{-1/2}^{1/2}
  \int_{\sqrt{1-x^{2}}}^{\infty}
  f(x+iy) \overline{g(x+iy)} y^{k} \, \frac{dx \, dy}{y^{2}}.
\]
It is well-known (see Theorem 6.12 of \cite{Iwa} for a proof) that
the Hecke operators are self-adjoint with respect to the Petersson inner
product, and this fact, together with the commutativity of $T_{p}$ and $T_{q}$,
implies that there is a basis for $S_{k}$ consisting of Hecke eigenforms,
each normalized so that the coefficient of $q$ is equal to $1$.

If
\[
  g(z) = \sum_{n=1}^{\infty} a(n) q^{n}
\]
is such a Hecke eigenform, Deligne proves in \cite{Del} that if $p$ is prime,
then
\[
  |a(p)| \leq 2 p^{\frac{k-1}{2}},
\]
as a consequence of the Weil conjectures. It follows from this that
$|a(n)| \leq d(n) n^{\frac{k-1}{2}}$ for all $n \geq 1$.

The self-adjoint property of the Petersson inner product implies
that if $g_{i}$ and $g_{j}$ are two distinct Hecke eigenforms, then $\langle
g_{i}, g_{j} \rangle = 0$. On the other hand, the second
equation on p. 251 of \cite{Iwa} gives that
\[
  L({\rm Sym}^{2} g_{i}, 1) = \frac{\pi^{2}}{6}
  \cdot \frac{(4 \pi)^{k} \langle g_{i}, g_{i} \rangle}{\Gamma(k)}.
\]
Here, $L({\rm Sym}^{2} g_{i}, s)$ is the symmetric square
$L$-function.  In the appendix to \cite{GHL}, Goldfeld, Hoffstein and
Lieman proved that $L({\rm Sym}^{2} g_{i}, s)$ has no Siegel zeroes,
and in \cite{Rou}, the second author used this to derive the lower
bound
\[
  L({\rm Sym}^{2} g_{i}, 1) \geq \frac{1}{64 \log(k)}.
\]

\section{Proof of Theorem~\ref{bound}}
\label{proofofbound}

Let $\ell = \dim S_{k}$ and write $k = 12 \ell + k'$, where $k' \in
\{ 0, 4, 6, 8, 10, 14 \}$. For each integer $m$ with $1 \leq m \leq
\ell$, we let $F_{k,m}(z)$ denote the unique weight $k$ modular form
with a Fourier expansion of the form
\[
  F_{k,m}(z) = q^m + \sum_{n=\ell+1}^{\infty} A_{k}(m,n) q^{n}.
\]

In \cite{DJ}, Duke and the first author gave a generating function for
the $F_{k,m}(z)$.  Note that the notation in this paper differs
slightly from theirs; $F_{k, m}$ is equal to the modular form $f_{k,
  -m}$ in~\cite{DJ}.

\begin{thmnonum}[Lemma~2 of \cite{DJ}]
We have
\[
  F_{k,m}(z) = \frac{1}{2 \pi i} \oint_{C} \frac{\Delta^{\ell}(z) E_{k'}(z)
  E_{14-k'}(\tau)}{\Delta^{1+\ell}(\tau) (j(\tau) - j(z))} p^{m-1} \,
  dp,
\]
where $p = e^{2 \pi i \tau}$ and $C$ denotes a (counterclockwise)
circle in the $p$-plane with sufficiently small radius.
\end{thmnonum}

Inspection of the integrand shows that the only poles of the
integrand (as $\tau$ varies) occur when $\tau$ is equivalent to $z$
under the action of $\SL_{2}(\Z)$. We change variables by setting
$\tau = u + iv$, $p = e^{2 \pi i \tau}$, $dp = 2 \pi i e^{2 \pi i
\tau}$, and let $v$ and $y$ be fixed constants. This gives
\[
  F_{k,m}(z) = \int_{-.5}^{.5} \frac{\Delta^{\ell}(z) E_{k'}(z) E_{14 - k'}(\tau)}
  {\Delta^{1+\ell}(\tau) (j(\tau) - j(z))} e^{2 \pi i m \tau} \, du,
\]
which is valid provided no point with imaginary part at least $v$
is equivalent to $z$ under the action of $\SL_{2}(\Z)$. It follows that
\[
  A_{k}(m,n) = \int_{-.5}^{.5} \int_{-.5}^{.5}
  \frac{\Delta^{\ell}(z) E_{k'}(z) E_{14 - k'}(\tau)}{\Delta^{1+\ell}(\tau)
  (j(\tau) - j(z))} e^{2 \pi i m \tau} e^{-2 \pi i n z} \, du \, dx,
\]
provided no point $\tau$ with $\IM{\tau} \geq v$ is equivalent to
any point $z$ with $\IM{z} = y$.

From this, it is clear that we can take absolute values to obtain
the bound \[\abs{A_k(m, n)} \leq \max_{\abs{u}, \abs{x} \leq .5}
\abs{\frac{\Delta(z)}{\Delta(\tau)}}^\ell
\abs{\frac{E_k'(z)E_{14-k'}(\tau)}{\Delta(\tau)(j(\tau)-j(z))}}
e^{-2\pi m v} e^{2\pi n y}.\]

Since $\Delta(z) = q - 24 q^2 + O(q^3)$, we have $\abs{\Delta(z)}
\leq e^{-2\pi y} + 24 e^{-4 \pi y} + B$, where $B$ is a bound on the
tail $\sum_{n=3}^\infty \tau(n) q^n$ of the series.  We can bound
the tail by $\sum_{n=3}^\infty d(n) n^{11/2} e^{-2\pi n y}$; using
the bound $d(n) \leq 2 \sqrt{n}$, we can exactly evaluate the sum
that results in terms of $y$.  This gives us an explicit upper bound
for $\abs{\Delta(z)}$ in terms of $y$.  Similarly, we find an lower
bound for $\abs{\Delta(\tau)}$ in terms of $v$.

For each of the six choices of $k'$, we bound $\abs{E_k'(z)
E_{14-k'}(\tau)}$ in terms of $y$ and $v$ by noting that
$\sigma_{k-1}(n) \leq 2\sqrt{n} n^{k-1} \leq 2n^{k}$, so that
\[\abs{E_k(z)} = \abs{1 - \frac{2k}{B_k} \sum_{n=1}^\infty
\sigma_{k-1}(n) q^n} \leq 1 + \frac{2k}{\abs{B_k}} \sum_{n=1}^\infty
2n^k e^{-2\pi n y}.\] This latter sum may be exactly evaluated in
terms of $y$.

At this point, we set $y = .865$ and $v = 1.16$; these values
satisfy the conditions above, since all points equivalent to $z = x
+ .865i$ under the action of $\SL_2(\Z)$ have imaginary part less
than 1.16, and give reasonable bounds for the quantities we are
studying. With these choices, we find that
\[\abs{\frac{\Delta(z)}{\Delta(\tau)}} \leq 7.358,\]
\[\abs{\frac{1}{\Delta(\tau)}} \leq 1488.802,\]
\[\abs{E_k'(z)E_{14-k'}(\tau)} \leq 40.368.\]

It remains to bound the quantity $\abs{j(\tau)-j(z)}$ on the
appropriate intervals.  We bound the tails of the two series, taking
all terms with exponent 10 and above for $j(z)$ and all terms with
exponent $5$ and above for $j(\tau)$.  Using the bounds given
in~\cite{BP}, we find that the tail of $j(z)$ is bounded by
\[\sum_{n=10}^\infty e^{-2\pi n (.865)} \frac{1}{\sqrt{2} n^{3/4}}
e^{4\pi \sqrt{n}} \left(1 - \frac{3}{32\pi\sqrt{n}} +
\frac{.055}{n}\right) \leq \frac{1.055}{\sqrt{2}} \sum_{n=10}^\infty
e^{-2\pi \sqrt{n} (.865\sqrt{n}-2)} \]\[\leq
\frac{1.055}{\sqrt{2}}\sum_{n=10}^\infty e^{-2\pi
\sqrt{n}(.2\sqrt{n})}  \leq .000003636545.\] Similarly, the tail of
$j(\tau)$ is bounded by $.000003636545$.

We now bound the main terms of $\abs{j(\tau)-j(z)}$.  Writing $j(z)
= q^{-1} + \sum c(n) q^n$, we must find a lower bound for
\[G(x, u) = \left| p^{-1} + \sum_{i=1}^4 c(i)p^i - q^{-1} -
\sum_{i=1}^9 c(i)q^i
\right|,\] where $p = e^{2\pi i (u + 1.16 i)}$, $q = e^{2 \pi i (x +
.865 i)}$, and $\abs{u}, \abs{x} \leq .5$.

To bound $G(x, u)$, we examine the function $G(x, u)^2$, which can
be written as an expression in $\cos(2\pi n x), \cos(2\pi n u),
\sin(2\pi n x)$, and $\sin(2\pi n u)$.  After finding bounds on the
partial derivatives of $G^2$ with respect to $x$ and $u$, we compute
its values on a grid of points satisfying $\abs{u}, \abs{x} \leq .5$
to see that $G^2 \geq 900$, implying that $G(x, u) \geq 30$ in this
range. The computations were performed using Maple, and were
shortened by noting that $G(x, u) = G(-x, -u)$; the bounds on
derivatives were calculated by trivially bounding the second
derivatives and, again, computing values on a grid of points.

Putting together these computations, we see that \[ |A_{k}(m,n)|
\leq 2003.34 \cdot 7.358^\ell e^{-2 \pi m \cdot 1.16}
  e^{2 \pi n \cdot 0.865}.\]

We now use this estimate on $|A_{k}(m,n)|$ to estimate $\langle G, G
\rangle$, where $G = \sum_{m=1}^{\ell} a(m) F_{k, m}$. We have
\begin{align*}
  \langle G, G \rangle &=
  \frac{3}{\pi} \int_{-1/2}^{1/2} \int_{\sqrt{1-x^{2}}}^{\infty}
  |G(x+iy)|^{2} y^{k-2} \, dy \, dx\\
  &\leq \frac{3}{\pi} \int_{\sqrt{3}/2}^{\infty}
  \int_{-1/2}^{1/2} |G(x+iy)|^{2} y^{k-2} \, dx \, dy.\\
\end{align*}
Plugging in the Fourier expansion $G(z) = \sum_{n=1}^{\infty} a(n)
q^{n}$ and using the fact that we are integrating over a complete
period gives
\[
  \langle G, G \rangle \leq
  \frac{3}{\pi} \sum_{n=1}^{\infty} |a(n)|^{2} \int_{\sqrt{3}/2}^{\infty}
  y^{k-2} e^{-4 \pi n y} \, dy.
\]
Setting $u = 4 \pi n y$, $du = 4 \pi n \, dy$ gives
\begin{equation}\label{ff}
  \langle G, G \rangle \leq
  \frac{12}{(4 \pi)^{k}} \sum_{n=1}^{\infty}
  \frac{|a(n)|^{2}}{n^{k-1}} \int_{2 \pi \sqrt{3} n}^{\infty} u^{k-2} e^{-u} \, du.
\end{equation}
We have
\[
  a(n) = \sum_{m=1}^{\ell} a(m) A_{k}(m,n)
\]
and so for $n \geq \ell + 1$, we have
\[
  |a(n)|^{2} \leq (2003.34)^{2} (7.358)^{2 \ell}
  \left| \sum_{m=1}^{\ell} a(m) e^{-2 \pi m \cdot 1.16}\right|^{2}
  \cdot e^{4 \pi n \cdot 0.865}.
\]

For $1 \leq n \leq \ell$ we use the simple bound
\[
  \int_{2 \pi \sqrt{3} n}^{\infty} u^{k-2} e^{-u} \, du
  \leq \int_{0}^{\infty} u^{k-2} e^{-u} \, du = (k-2)!.
\]
Hence, the contribution to $\langle G, G \rangle$ from the terms
with $1 \leq n \leq \ell$ is at most
\[
  \frac{12 (k-2)!}{(4 \pi)^{k}} \sum_{n=1}^{\ell} \frac{|a(n)|^{2}}{n^{k-1}}.
\]

For $n \geq \ell + 1$ we use that
\[
  \int_{2 \pi \sqrt{3} n}^{\infty} u^{k-2} e^{-u} \, du
  = e^{-2 \pi \sqrt{3} n}
\sum_{i=0}^{k-2} \frac{(k-2)!}{i!} (2 \pi \sqrt{3} n)^{i}.
\]
Since the highest power of $n$ in this expression is $k-2$, the
piece
\[
  \frac{1}{n^{k-1}} \sum_{i=0}^{k-2} \frac{(k-2)!}{i!} (2 \pi \sqrt{3} n)^{i}
\]
of the right side of equation~\eqref{ff} is a decreasing function of
$n$ and is therefore bounded by
\[
  \frac{1}{(\ell+1)^{k-1}} \sum_{i=0}^{\infty} \frac{(k-2)!}{i!}
  (2 \pi \sqrt{3} (\ell+1))^{i}
  = \frac{(k-2)! e^{2 \pi \sqrt{3} (\ell + 1)}}{(\ell+1)^{k-1}}.
\]
Hence, the contribution to $\langle G, G \rangle$ from the terms with
$n \geq \ell + 1$ is at most
\[
  \frac{12}{(4 \pi)^{k}} \cdot (2003.34)^{2} (7.358)^{2 \ell}
  \left| \sum_{m=1}^{\ell} a(m) e^{-2 \pi m \cdot 1.16} \right|^{2}
  \cdot \frac{(k-2)! e^{2 \pi \sqrt{3} (\ell + 1)}}{(\ell+1)^{k-1}}
  \cdot \sum_{n=\ell+1}^{\infty} e^{4 \pi n \cdot 0.865}
  e^{-2 \pi \sqrt{3} n}.
\]
The sum on $n$ is a geometric series, and we have $4 \pi \cdot 0.865
- 2 \pi \sqrt{3} \leq -0.01288$. This gives the bound
\[
  \frac{(k-2)! (12168805)^{2}}{(4 \pi)^{k}}
  \left| \sum_{m=1}^{\ell} a(m) e^{-2 \pi m \cdot 1.16} \right|^{2}
  \cdot \frac{(7.358)^{k/6} 12^{k} e^{k \pi \sqrt{3}/6}
  e^{-0.00107 k}}{k^{k-1}}.
\]

Thus, we have
\[
  \langle G, G \rangle \leq
  \frac{12 (k-2)!}{(4 \pi)^{k}}
  \sum_{m=1}^{\ell} \frac{|a(m)|^{2}}{m^{k-1}}
  + \frac{(12168805)^{2} (k-2)!}{(4 \pi)^{k}}
  \left| \sum_{m=1}^{\ell} a(m) e^{-2 \pi m \cdot 1.16}\right|^{2}
  \cdot \frac{(41.41)^{k}}{k^{k-1}}.
\]
Now, we write $G = \sum_{i=1}^{\ell} c_{i} g_{i}$, where the $g_{i}$
are the normalized Hecke eigenforms. Using the lower bound on
$L({\rm Sym}^{2} g_{i}, 1)$ and the relation between $L({\rm Sym}^{2} g_{i}, 1)$ and $\langle g_{i}, g_{i}
\rangle$, we get
\begin{align*}
  \langle G, G \rangle &= \sum_{i=1}^{\ell} |c_{i}|^{2} \langle g_{i}, g_{i}
  \rangle\\
  &\geq \sum_{i=1}^{\ell} |c_{i}|^{2} \cdot \left(\frac{3 (k-1)!}{32 \pi^{2} (4 \pi)^{k} \log(k)}\right).\\
\end{align*}
This gives an upper bound on $\sum_{i=1}^{\ell} |c_{i}|^{2}$ in terms of
$\langle G, G \rangle$. The Cauchy-Schwarz inequality gives
\begin{align*}
  \sum_{i=1}^{\ell} |c_{i}| \leq
  & \sqrt{\ell} \sqrt{\sum_{i=1}^{\ell} |c_{i}|^{2}}\\
  & \leq \sqrt{\frac{k}{k-1} \cdot \frac{32 \pi^{2}}{3} \log(k)}
  \cdot \sqrt{\sum_{m=1}^{\ell} \frac{|a(n)|^{2}}{m^{k-1}}}\\
  &+ \sqrt{\frac{k}{k-1} \cdot \frac{32 \pi^{2}}{3} \cdot \log(k)}
  \cdot 12168805 \cdot \left|\sum_{m=1}^{\ell} a(m) e^{-7.288m}\right|
  \cdot \frac{(41.41)^{k/2}}{k^{(k-1)/2}}\\
  &\leq \sqrt{\log(k)} \left(11 \cdot \sqrt{\sum_{m=1}^{\ell}
  \frac{|a(m)|^{2}}{m^{k-1}}}
  + \frac{e^{18.72} (41.41)^{k/2}}{k^{(k-1)/2}}
  \cdot \left|\sum_{m=1}^{\ell} a(m) e^{-7.288m} \right|\right).
\end{align*}
This concludes the proof of Theorem~\ref{bound}.

\section{Proof of Theorem~\ref{extreme}}
\label{proofofextreme}

Write $F_{k,0}(z) = E_{k}(z) + h(z)$, where
\[
  h(z) = \sum_{n=1}^{\infty} b(n) q^{n}.
\]
Since $F_{k,0}(z) = 1 + O(q^{\ell + 1})$, we have
\[
  b(m) = \frac{2k}{B_{k}} \sigma_{k-1}(m)
\]
for $1 \leq m \leq \ell$. We now apply Theorem~\ref{bound}, which gives
that $b(n)$ is bounded by
\[
  \sqrt{\log(k)} \left(11
  \sqrt{ \sum_{m=1}^{\ell} \frac{|b(m)|^{2}}{m^{k-1}} }
  + \frac{e^{18.72} (41.41)^{k/2}}{k^{(k-1)/2}}
  \cdot \left| \sum_{m=1}^{\ell} b(m) e^{-7.288 m} \right|\right) d(n) n^{\frac{k-1}{2}}.
\]
We have that
\[
  \zeta(k) = \frac{(-1)^{\frac{k}{2} - 1} (2 \pi)^{k} B_{k}}{(k-1)! \cdot 2k}.
\]
If $k \geq 12$, then $1 \leq \zeta(k) \leq \zeta(12) \leq 1.00025$.
Thus, for $k \geq 12$ we have
\[
  0.9997 \frac{(2 \pi)^{k}}{(k-1)!}
\leq -\frac{2k}{B_{k}} \leq \frac{(2 \pi)^{k}}{(k-1)!}.
\]
Now, we have
\[
  \sigma_{k-1}(m) = \sum_{d | m} d^{k-1} =
  \sum_{d | m} (m/d)^{k-1} = m^{k-1} \sum_{d | m} \frac{1}{d^{k-1}}
  \leq m^{k-1} \zeta(k-1).
\]
We have
\[
  \sqrt{\sum_{m=1}^{\ell} \frac{|b(m)|^{2}}{m^{k-1}}}
  \leq -\frac{2k \zeta(k-1)}{B_{k}} \sqrt{\sum_{m=1}^{\ell} m^{k-1}}.
\]
Also,
\begin{align*}
  \sum_{m=1}^{\ell} m^{k-1} &= \int_{1}^{\ell + 1} \lfloor x \rfloor^{k-1} \, dx
  \leq \int_{1}^{\ell + 1} x^{k-1} \, dx \leq \frac{(\ell+1)^{k}}{k}\\
  &\leq \frac{\ell^{k} \left(1 + \frac{1}{\ell}\right)^{12 \ell + 12}}{k}
  \leq \frac{e^{12} \ell^{k}}{k} \cdot \left(1 + \frac{1}{\ell}\right)^{12}.
\end{align*}
Thus, the contribution from the first
term in Theorem~\ref{bound} is
\[
  \frac{(2 \pi)^{k}}{(k-1)!} \frac{11 \cdot 1.0005 \cdot e^{6} \ell^{k/2}}{\sqrt{k}}
  \left(1 + \frac{1}{\ell}\right)^{6} \sqrt{\log(k)}.
\]

The function $m^{k-1} e^{-7.288 m}$ always has a maximum at $m
= \ell$. Thus, the second term of the bound from
Theorem~\ref{bound} is at most
\begin{align*}
  & \frac{\zeta(11) e^{18.72} (41.41)^{k/2} \ell^{k} e^{-7.288 \ell}
  (2 \pi)^{k} \sqrt{\log(k)}}{(k-1)! k^{(k-1)/2}}\\
  &\leq \frac{(2 \pi)^{k}}{(k-1)!} e^{28.4657} \ell^{(k+1)/2} (1.0242382)^{k/2}
  \sqrt{\log(k)}.
\end{align*}

Adding the two contributions above, we have that
\[
  C(h) \leq \frac{(2 \pi)^{k}}{(k-1)!} e^{28.466} \sqrt{\ell \log(k)}
  (1.0242382 \ell)^{k/2},
\]
and so $|b(n)| \leq C(h) d(n) n^{\frac{k-1}{2}} \leq
2C(h) n^{k/2}$. Now, we have
\[
  a(n) = -\frac{2k}{B_{k}} \sigma_{k-1}(n)
  + b(n) \geq 0.9997 \frac{(2 \pi)^{k}}{(k-1)!}
  n^{k-1} - 2 C(h) n^{k/2}.
\]
The right hand side is positive if
\begin{align*}
  n^{\frac{k}{2} - 1} &\geq \frac{2 e^{28.466} \sqrt{\ell \log(k)} (1.0242382 \ell)^{k/2}}{0.9997}\\
  n &\geq e^{58.366/(k-2)} \left(\ell^{3} \log(k)\right)^{\frac{1}{k-2}}
  \cdot 1.0242382 \ell.
\end{align*}
This concludes the proof of Theorem~\ref{extreme}.

\section{Proof of Corollary~\ref{cor0}}
\label{computations}

To verify that all Fourier coefficients of $F_{k,0}(z)$ are
non-negative for\\ $k \in \{81288, 81460, 81632\}$, we use the bound
from Theorem~\ref{extreme}. This shows that any negative Fourier
coefficient occurs within the first 10000. We find the unique linear
combination
\[
  \sum_{i=0}^{k/4} c_{i} E_{4}^{k - 3i} \Delta^{i} = 1 + O(q^{\ell + 1})
\]
and this form will equal $F_{k,0}(z)$. It then suffices to check the
first 10000 Fourier coefficients are non-negative. These
computations are performed in Magma \cite{Magma}, and take
approximately 3 days for each weight.

Recall that
\[
  F_{k,0}(z) = \sum_{n=0}^{\infty} a(n) q^{n}.
\]
We will show that $a(\ell + 2) < 0$ for $k$ sufficiently large
(depending on $k \bmod 12$), making effective the work of Mallows,
Odlyzko, and Sloane. Write
\[
  E_{4}^{-k/4} = \sum_{n=0}^{\infty} A(n) j^{-n}
\]
where $j$ is the usual $j$-function. B\"urmann's theorem gives that
\begin{equation}
\label{readoffcoeff}
  A(n) = \left(-\frac{k}{4n}\right) \cdot \text{ the coefficient of }
  q^{n-1} \text{ in } \left(\frac{dE_{4}}{dq} \frac{E_{4}^{3n - k/4 - 1} q^{n}}{\Delta^{n}}\right).
\end{equation}
Mallows, Odlyzko, and Sloane show (see \cite{MOS}, pg. 73) that
\begin{align*}
  a(\ell + 1) &= -A(\ell + 1) > 0\\
  a(\ell + 2) &= -A(\ell + 2) + A(\ell + 1) \left(24 \ell - 240 \nu + 744\right).
\end{align*}
We write
\begin{align*}
  A(\ell + 1) &= -\frac{k}{4(\ell + 1)} \int_{-1/2}^{1/2} \theta(E_{4}) E_{4}^{2 - \nu}
  \frac{1}{\Delta^{\ell + 1}} \, dx\\
  A(\ell + 2) &= -\frac{k}{4(\ell + 2)} \int_{-1/2}^{1/2} \theta(E_{4}) E_{4}^{5 - \nu} \frac{1}{\Delta^{\ell + 2}} \, dx
\end{align*}
where $\theta\left(\sum a_{n} q^{n}\right) = \sum n a_{n} q^{n}$,
and the integrals are over the line segment $x + iy$, $-1/2 \leq x
\leq 1/2$ where $y$ is fixed. We wish to find an upper bound on
$|A(\ell + 2)|$ and a lower bound on $|A(\ell + 1)|$.

We choose $y$ so that $\frac{\Delta'(iy)}{\Delta(iy)} = 0$ (so $y
\approx 0.52352$). We write the integrals above in the form
\[
  \int_{-1/2}^{1/2} H_{j}(x+iy) e^{-(\ell + j) \ln(\Delta(x+iy))} \, dx
\]
where $H_{1}(x+iy) = \theta(E_{4})(x + iy) E_{4}(x+iy)^{2 - \nu}$
and $H_{2}(x+iy) = \theta(E_{4})(x+iy) E_{4}(x+iy)^{5-\nu}$.

If $B(x) = -\ln(\Delta(x+iy))$, then $|B(x)| \leq B(0) \approx
4.23579$. Moreover, the choice of $y$ gives that $B'(0) = 0$. We use
Taylor's theorem with the Lagrange form of the remainder to write
\[
  B(x) = B(0) + \frac{1}{2} x^{2} {\rm Re}(B)''(z_{1}) + \frac{i}{2} x^{2} {\rm Im}(B)''(z_{2}) := B(0) + x^2 C_{1}(x) + i x^{2} C_{2}(x).
\]
for some $z_{1}$ and $z_{2}$ between $0$ and $x$. We bound from
above and below the second derivatives of the real and imaginary
parts of $B$. We derive similar bounds on $H_{1}(x+iy)$ and
$H_{2}(x+iy)$.

We then have
\[
  e^{-(\ell + j) B(x)} = e^{-(\ell + j) B(0)} \cdot e^{C_{1}(x) x^{2}} \left(\cos\left((\ell + j) C_{2}(x) x^{2}\right) + i \sin\left((\ell + j) C_{2}(x) x^{2}\right)\right).
\]
Since the integrals we are studying are both real, we wish to
approximate the real part of the integrand. The main contribution
comes in an interval of length about $\frac{1}{\sqrt{\ell}}$ in a
neighborhood of $x = 0$, chosen so that $\cos((\ell + j) C_{2}(x)
x^{2})$ is positive. We bound the contribution of the remaining part
of $-1/2 \leq x \leq 1/2$ trivially.

The bounds we obtain from this method show that $a(\ell + 2) < 0$ if
$k \geq 84636$, $k \geq 83332$, and $k \geq 82532$ if $\nu = 0$,
$\nu = 1$, or $\nu = 2$, respectively. We use \eqref{readoffcoeff}
to compute the coefficient $a(\ell + 2)$ for all $k$ between the
bounds given in Corollary~\ref{cor0} and the bounds above. This
concludes the proof.

\providecommand{\bysame}{\leavevmode\hbox to3em{\hrulefill}\thinspace}
\providecommand{\MR}{\relax\ifhmode\unskip\space\fi MR }
\providecommand{\MRhref}[2]{%
  \href{http://www.ams.org/mathscinet-getitem?mr=#1}{#2}
}
\providecommand{\href}[2]{#2}

\end{document}